\documentstyle[12pt]{article}

  \newcommand{\const}{\rm const}
  
  \newcommand{\supp}{\rm supp}

  \newcommand{\Dom}{\rm  Dom}

  \newcommand{\Sub}{\rm Sub}

\begin{document}

   \begin{center}

{\bf  Generalization of tail inequalities for random variables,}\\

\vspace{4mm}

{\bf using in  the martingale theory. }\\

\vspace{5mm}

{\bf M.R.Formica, \ E.Ostrovsky, \ L.Sirota. } \\

 \end{center}

 \vspace{5mm}

\ Universit\`{a} degli Studi di Napoli Parthenope, via Generale Parisi 13, Palazzo Pacanowsky, 80132,
Napoli, Italy. \\
e-mail: \ mara.formica@uniparthenope.it \\

 \ Department of Mathematics and Statistics, Bar-Ilan University,\\
59200, Ramat Gan, Israel. \\
e-mail: \ eugostrovsky@list.ru\\

 \ Department of Mathematics and Statistics, Bar-Ilan University,\\
59200, Ramat Gan, Israel. \\
e-mail: \ sirota3@bezeqint.net \\

\vspace{5mm}

\begin{center}

 \ {\bf Abstract.}

\end{center}

 \hspace{3mm} We generalize a famous tail Doob's inequality, relative two non - negative random variables,
 arising in the martingale theory, in two directions:  on the more general source data and on the
 random variables belonging to the so - called Grand Lebesgue Spaces.\par
  \ We bring also several examples  in each sections in order to show the exactness of our estimates. \par

\vspace{5mm}

 \ {\it Key words and phrases.} Martingale, random variable (r.v.), expectation, generating function, inequalities,
 upper and lower estimates, moment, dilation operator, H\"older's inequality, event,
 tail of distribution, Lebesgue - Riesz. Grand Lebesgue and Orlicz spaces and norms, Young - Fenchel transform,
Doob's and other inequalities, probability space, examples.\\

\vspace{5mm}

\section{Notations. Statement of problem.}

\vspace{5mm}

 \hspace{3mm} Let $ \ (\Omega = \{\omega\}, \cal{B}, {\bf P}) \ $ be non - trivial probability space with expectation $ \ {\bf E}. \ $
 The classical Doob's inequality, see e.g. \cite{Doob 1}, \  \cite{Doob 2}, see also \cite{Pollard}, tell us that if the  non - negative 
 numerical valued r.v. - s $ \   (\xi, X) \ $ are such that for some positive finite constants $ \  \beta, C  \ $

\begin{equation} \label{source ineq}
t \ {\bf P} (\xi > \beta t) \le C \ {\bf E} [ X \ I(\xi > t) ], \ t \ge 0,
\end{equation}
here and further $ \ I(A) \ $ denotes as ordinary  an indicator function of the random event (predicate) $ \ A, \ $ then

\begin{equation} \label{classical ineq}
||\xi||_p =  [{\bf E} | \xi^p|]^{1/p} \le C \ \frac{p}{p - 1} \ \beta^p \ |X|_p, \ p > 1.
\end{equation}

 \ Henceforth $ \ ||\eta||_p \ $ denotes as ordinary the classical Lebesgue - Riesz norm of the r.v. $ \ \eta: \ $

$$
||\eta||_p := \left[ \ {\bf E}|\eta|^p \ \right]^{1/p}, \ p \ge 1.
$$

 \vspace{3mm}

 \ This inequality  (\ref{classical ineq}) play a very important role in particular, in the martingale theory, see
\cite{Burkholder 1}, \cite{Burkholder 2}. \par

\vspace{4mm}

 \ {\bf  Our aim in this short report is generalization of Doob's inequality in two  mentioned directions.  } \par

\vspace{4mm}

\begin{center}

 \ {\sc Statement of problem.} \par

\end{center}

\vspace{4mm}

 \hspace{3mm} Given:

\begin{equation} \label{modern ineq}
h(t) \ {\bf P} (\xi > \beta t) \le C \ {\bf E} [ X \ I(\xi > t) ], \ t \ge 0,
\end{equation}
 where $ \ h = h(t)   \ $ is certain non - negative continuous strictly increasing  deterministic function and as before  $ \ X \ $  be
non - negative random variable,  \ $ \ \beta, \ C = \const \in (0,\infty).\ $  \par

\vspace{3mm}

 \hspace{3mm} Let  also $ \  g = g(t) \ $ be some non - negative continuous differentiable deterministic strictly
increasing function such that $ \  g(0) = g(0+) = 0. \ $ We intent on the assumption of (\ref{modern ineq}) to find such a
function $ \ g = g(t) \ $ for which the following moment is finite:

\begin{equation} \label{target}
{\bf E} g(\xi) \le Z(\beta,h, X) < \infty
\end{equation}
and moreover to estimate the right - hand side of the relation (\ref{target}), i.e. the using in practice  functional $ \ Z = Z(\beta, h,X). \ $
We  bring also the examples in order to show the exactness of obtained estimates.\par

\vspace{5mm}

\section{Main result.}

\vspace{5mm}

 \hspace{3mm} We have

$$
{\bf P} (\xi/\beta > t) \le \frac{C}{h(t)} \ {\bf E} \{ \  X \ I(\xi > t) \ \}, \ t > 0,
$$
where as ordinary $ \ I(A) \ $ denotes the indicator function for the random event  (predicat) $ \ A. \ $  Hence

\begin{equation} \label{integr relation}
\int_0^{\infty} p t^{p-1} \ {\bf P} (\xi/\beta > t) \ dt \le \int_0^{\infty} \frac{ C \ p \ t^{p-1}}{h(t)} \ {\bf E} \{ \ X \ I(\xi > t)  \ \} \ dt, \ p \ge 1.
\end{equation}

 \ The left - hand side of (\ref{integr relation})  is equal to $ \  L := {\bf E}|\xi/\beta|^p = ||\xi||^p_p \ \beta^{-p}. \ $ Let us investigate the right - hand side $ \ R \ $
 of (\ref{integr relation}). We deduce by virtue of Fubini's theorem

$$
R = C \ p \ {\bf E} \ \left[ \ X \ \int_0^{\xi} \frac{t^{p-1}}{h(t)} \ dt  \ \right].
$$
 \hspace{3mm}  Let us introduce the new important random variable (measurable function), if  it there exists:

$$
\kappa_p(\xi) =  \kappa_p \ \stackrel{def}{=} \int_0^{\xi}  \ \frac{t^{p-1}}{h(t)} \ dt,
$$
then

$$
||\xi||^p_p \ \beta^{-p} \le C \ {\bf E} \ [X \ \kappa_p(\xi)].
$$
 \  Denote now

$$
K_p(\theta) := ||\kappa_p(\xi)||_{\theta}  < \infty, \   \exists \ \theta > 1; \ \alpha  = \alpha(\theta):= \theta/(\theta - 1)  = \theta' \ \in (1,\infty).
$$
and assume

\vspace{3mm}

\begin{equation} \label{K condition}
\exists \ \theta > 1 \ \Rightarrow \   K_p(\theta) := ||\kappa_p(\xi)||_{\theta}  < \infty,
\end{equation}
and denote by $ \ \Theta = \Theta(p) = \Theta(p)[h,\xi] \ $ the  set all the values $ \ \theta \ $ for which the value  $ \ K_p(\theta) \ $ is finite:

$$
\Theta = \Theta(p) = \Theta(p)[h,\xi]  = \{ \ \theta, \ \theta \ge 1, \ K_p(\theta) < \infty \ \}.
$$

\vspace{3mm}

 \ We apply as expected  the H\"older's inequality for such a values of the auxiliary  parameter $ \ \theta \in \Theta \ $

$$
||\xi||_p^p \le C \ \beta^p \ \ K_p(\theta) \ ||X||_{\alpha(\theta)}.
$$

\vspace{3mm}

 \hspace{3mm} {\it We impose  also the following condition on the source datum}

 \begin{equation} \label{key restriction}
 K_p(\theta) =  ||\kappa_p(\xi)||_{\theta} \le v(\theta,p,r) \ ||\xi||_p^r, \
 \end{equation}

 \begin{equation} \label{restric nonat}
 \exists \ r = r(\theta,p) \in [1, p), \ \exists \ v = \ v(\theta,p,r)  < \infty.
 \end{equation}

\vspace{3mm}

 \ Denote by $ \ R = R(\theta,p) \ $ the set of finiteness of the value $ \ r \ $ in the relation (\ref{key restriction}):

$$
R = R(\theta,p) := \{ r: \ 1 < r < p, \  K_p(\theta)  < \infty \}.
$$

\vspace{3mm}

 \hspace{3mm} To summarize: \par

 \vspace{4mm}

 \ {\bf Theorem 2.1.} We deduce under formulated above conditions

\vspace{3mm}

\begin{equation} \label{key estim}
||\xi||_p \le \left[ \ C \ v(\theta,p) \ \beta^p  \right]^{1/(p - r)} \cdot ||X||_{\alpha(\theta)}^{1/(p - r)}.
\end{equation}

\vspace{3mm}

 \ Of course,

\begin{equation} \label{key  opt estim}
||\xi||_p \le  \inf_{r \in R} \ \inf_{\theta \in \Theta} \  \left\{ \ C \ v(\theta,p,r) \ \beta^p\right\}^{1/(p - r)}  \cdot  ||X||_{\alpha(\theta)}^{1 / (p - r)}.
\end{equation}

\vspace{4mm}

\begin{center}

 {\sc Some examples.} \\

\end{center}

\vspace{4mm}

 \ {\bf Example 2.1.} Put now $ \  h(t) = t, \ $ so that

$$
t \ {\bf P}(\xi > \beta \ t) \le C {\bf E} (X \ I(\xi > t)), \ t, \xi, X \ge 0, \ p > 1.
$$

 \ We choose in the conditions of theorem 2.1  $ \ p =1, \ r = p/(p-1),  \ $ then

 $$
 v(t) = [p/(p-1) ] \cdot t^{p-1}.
 $$

 \ We get using the proposition of theorem 2.1,  after simple calculations,  the classical result, see e.g. \cite{Pollard}

\begin{equation} \label{beta}
|\xi|_p \le C \ \frac{p}{p-1} \ \beta^p \ |X|_p, \ p > 1.
\end{equation}

\vspace{4mm}

\ {\bf Example 2.2.}  A more general case. Suppose as before that $ \ \xi, \ X \ $ are non - negative r.v. such that

\vspace{3mm}

$$
\exists \ C = \const > 0, \hspace{3mm} \exists \Delta = \const > 1,
$$
and

\begin{equation} \label{cond Delta}
t^{\Delta} \ {\bf P}(\xi > \beta \ t) \le C {\bf E} (X \ I(\xi > t)), \ t \ge 0, \ p > \Delta.
\end{equation}

 \ We get again using the proposition of theorem 2.1  after simple calculations  the following estimation

\begin{equation} \label{beta Delta}
|\xi|_p \le C \ \frac{p}{ p - \Delta} \ \beta^p \ |X|_p, \ p > \Delta.
\end{equation}

\vspace{5mm}

\section{Unimprovability of our estimations. Lower bounds.}

\vspace{5mm}

 \hspace{3mm} Let us show the exactness of obtained results, in particular, ones \ (\ref{beta}), \ (\ref{beta Delta}). Introduce
 the following important functionals,  which are  responsible for the lower estimate.

$$
Y[C,\beta, \Delta](\xi, X,p) \stackrel{def}{=} \left[ \ \frac{|\xi|_p}{ C \ p \ (p - \Delta)^{-1} \ \beta^p \ |X|_p} \ \right],
$$

\begin{equation} \label{K functional}
 U = U[C, \beta,\Delta,p]  =  \stackrel{def}{=} \sup_{p > \Delta} \ \sup_{\xi \in L_p} \ \sup_{X \in L_p} Y[C,\beta, \Delta](\xi, X,p),
\end{equation}
where  all the supremums are calculated over the r.v. - s $ \ \xi, \ X \ $
 satisfying the condition (\ref{cond Delta}) and when $ \ p > \Delta, \ \Delta = \const > 1. \ $ \par

\vspace{4mm}

 \ {\bf Proposition 3.1.} \\

 \vspace{3mm}

\begin{equation} \label{lower bound}
U(1,1, \Delta) = 1.
\end{equation}

 \vspace{3mm}

 \hspace{3mm} {\bf Proof.} The upper estimate  $ \  U(C,\beta, \Delta) \le 1 \ $ is contained in (\ref{beta Delta}). In order to deduce the lower one,
 we bring an example. \par
 \ Let us choose
$ \ C, \beta = 1 \ $ and bring as the variables $ \ \xi, \ X \ $ the following: $ \ X_0 = \xi_0 \ $ and let the random  positive variable $ \ \xi_0, \ $
as well as one $ \ X_0, \ $ has a standard  exponential distribution

$$
{\bf P} (\xi_0 > t) = {\bf P}(X_0 > t) = e^{-t}, \  \ \xi_0, \ X_0 > 0, \ t > 0;
$$
then the natural generating function for these r.v.- s has a form

$$
\nu(p) \stackrel{def}{=}  |\xi_0|_p = |X_0|_p = \Gamma^{1/p} (p + 1), \ p \ge 1,
$$
where as ordinary $ \ \Gamma(\cdot) \ $ is Euler's Gamma function.\par

\vspace{3mm}

\ Note that as $ \ p \to \infty \ \Rightarrow \nu(p) \sim p/e. $ \  Note also that the condition (\ref{source ineq})
is satisfied.\par

\vspace{3mm}

 \ We have

$$
Y[C,\beta, \Delta](\xi_0, X_0,p) = \frac{p- \Delta}{p}, \ p > \Delta.
$$

 \ Following,

\begin{equation} \label{frac xi X}
U(1,1,\Delta) \ge \sup_{p > \Delta} \left\{ \ \frac{p - \Delta}{p} \ \right\}.
\end{equation}

 \ Our proposition (\ref{lower bound}) follows immediately from the equality

$$
\lim_{p \to \infty} \left\{ \ \frac{p - \Delta}{p} \ \right\} = 1.
$$

\vspace{3mm}

 \ {\bf Remark 3.1.} \ The cases $ \ C,\beta \ne 1 \ $  may be considered quite analogously. \par

\vspace{5mm}

\section{ Generalization on the Grand Lebesgue Spaces approach. Examples.}

\vspace{5mm}

 \hspace{3mm} We intent in this section to extend the previous results upon the so - called Grand Lebesgue Spaces  (GLS) of the random variables. \par

\vspace{4mm}

\hspace{3mm} Let $ \  (a,b) = \const, 1 \le a < b \le \infty.  \  $ Let also $ \ \psi = \psi(p), \ p \in (a,b) \ $  be certain numerical valued
 measurable strictly positive: $ \  \inf_{p \in (a,b)} \psi(p) > 0   \ $  function, not necessary to be bounded. Denotation:
 $ \ \Dom (\psi) \stackrel{def}{=} \{ \ p: \ \psi(p) < \infty \ \}, \ $

 $$
  (a,b) := \supp (\psi); \  \Psi(a,b) := \{ \  \psi: \ \supp (\psi) = (a,b)  \ \},
 $$

$$
\Psi \stackrel{def}{=} \cup_{(a,b) } \Psi(a,b).
$$

\vspace{3mm}

 \hspace{3mm} {\bf Definition 4.1.}, see e.g.  \cite{KosOs}, \cite{Ermakov}, \cite{Fiorenza-Formica-Gogatishvili-DEA2018}.
 \ Let the function $ \  \psi = \psi(p), \ p \in (a,b) \ $ belongs to the set
 $ \   \Psi(a,b): \hspace{3mm} \psi(\cdot) \in \Psi(a,b), \ $  which is named as {\it generating function} for introduced after space.
 \ The  so - called {\it Grand Lebesgue Space} $ \  G \psi \ $
 is defined as a set of all random variables (measurable functions)  $ \ \tau \ $ having a finite norm

\vspace{3mm}

\begin{equation} \label{def GLS norm}
||\tau||G\psi  \stackrel{def}{=}  \sup_{p \in (a,b)} \left\{ \ \frac{||\tau||L_p(\Omega)}{\psi(p)} \ \right\} =  \sup_{p \in (a,b)} \left\{ \ \frac{||\tau||_p}{\psi(p)} \ \right\}.
\end{equation}

\vspace{4mm}

 \ The particular case of these spaces  and under  some additional restrictions on the generating function $ \ \psi = \psi(p) \ $
  correspondent to the so - called {\it Yudovich spaces,}  see  \cite{Yudovich 1}, \cite{Yudovich 2}.
These spaces was applied at first in the theory of Partial Differential Equations (PDE),  see  \cite{Chen},  \cite{Crippa}. \par

\vspace{3mm}

 \  These spaces are  complete Banach functional rearrangement invariant; they are
 investigated in many works, see e.g. \cite{ErOs}, \cite{Ermakov}, \cite{Fiorenza2}, \cite{Fiorenza-Formica-Gogatishvili-DEA2018},
\cite{fioforgogakoparakoNAtoappear}, \cite{fioformicarakodie2017}, \cite{formicagiovamjom2015}, \cite{Formica Ostrovsky Sirota weak dep},
\cite{KosOs}, \cite{KozOsSir2017}, \cite{KosOs equivalence}, \cite{Liflyand}, \cite{Ostrovsky1}. It is important for us in particular to note
that there is exact of course up to finite multiplicative constant
interrelations  under certain natural conditions on the generating function between belonging the r.v. $ \ \tau \ $ to this
space and it tail behavior. Indeed, assume for the definiteness that $ \  \tau \in G \psi  \ $ and moreover $ \ ||\tau||G\psi = 1;  \ $ then

\begin{equation} \label{tail estim}
T_{\tau}(t) \le  \exp \{ \  -  h^*(\ln t)  \  \}, \ t \ge e,
\end{equation}
where $ \ h(p) = h[\psi](p) := p \ln \psi(p) \ $  and $ \  h^*(\cdot)   \ $  is the famous Young - Fenchel (Legendre) transform of the function $ \ h(\cdot): \ $

$$
h^*(u) \stackrel{def}{=} \sup_{p \in \Dom( \psi)} (pu - h(p)).
$$

 \ Inversely, let the tail function $ \ T_{\tau}(t), \ t \ge 0  \ $  be given.  Introduce the following so - called {\it natural function}
 generated by $ \ \tau \ $

\begin{equation} \label{natural function}
 \psi_{\tau}(p) \stackrel{def}{=} \left[ \   p \int_0^{\infty}  \ t^{p-1} \ T_{\tau}(t)  \ dt \ \right]^{1/p} = ||\tau||L_p(\Omega),
\end{equation}
 if it is finite for some value $ \ b \in (a,\infty], \ $ following, it is finite at last for all the values $ \ p \in (a,b). \ $ \par

 \ As long as

$$
{\bf E} |\tau|^p = p \int_0^{\infty}  \ t^{p-1} \ T_{\tau}(t)  \ dt \ = \psi^p_{\tau}(p), \ p \in [1,b),
$$
 we conclude that if the last {\it natural} for the r.v. $ \ \tau \ $ function $ \ \psi_{\tau}(p) \ $ is finite inside some non - trivial
 segment $ \ p \in [1,b), \ 1 < b \le \infty, \ $   then

$$
\tau \in G\psi_{\tau}; \ \hspace{3mm} ||\tau||G\psi_{\tau} = 1.
$$

\vspace{5mm}

 \ Further, let the estimate (\ref{tail estim}) be given. Suppose in addition   that the generating function $ \ \psi = \psi(p), \ p \in \Dom(\psi) \ $
 is continuous and suppose in the case when $ \ b = \infty \ $

\begin{equation} \label{lim 0}
\lim_{p \to \infty} \frac{\psi(p)}{p} = 0.
\end{equation}

 \ Then the r.v. $ \ \tau \ $ belongs to the Grand Lebesgue Space $ \ G \psi: \    \ $

\begin{equation} \label{inverse est}
||\tau||G\psi \le K[\psi] < \infty,
\end{equation}
see e.g.   \cite{KosOs equivalence}.\par

 \ These conditions on the generating function $ \ \psi(\cdot) \ $  are satisfied for example for the functions  $ \ \psi_{m,L}(p) \ $ of the form

\begin{equation} \label{psi mL}
\psi_{m,L}(p) \stackrel{def}{=} p^{1/m} \ L(p), \ m = \const > 1, \ b = \infty,
\end{equation}
 where $ \ L = L(p) \ $ be some continuous strictly positive {\it  slowly varying }  at infinity function such that

\begin{equation} \label{theta condition}
\forall \theta > 0 \ \Rightarrow \sup_{p \ge 1} \left[  \ \frac{L(p^{\theta})}{L(p)} \ \right] = C(\theta) < \infty.
\end{equation}

 \ For instance, $ \ L(p) = [\ln(p+1) ]^r, \ r \in R. \ $ \par

\vspace{3mm}

 \ We conclude that under formulated restrictions the r.v. $ \ \tau \ $ belongs to the space $ \ G\psi_{m,L}: \ $

\begin{equation} \label{Gpsi m L}
\sup_{p \ge 1}  \left\{ \ \frac{||\tau||_{p,\Omega}}{\psi_{m,L}(p)}  \ \right\} = C(m,L) < \infty
\end{equation}

\ if and only if

\begin{equation} \label{tail behav}
T_{\tau}(u) \le \exp \left( \  - C_2(m,L) \ u^m /L(u) \ \right), \ u \ge e, \  \exists \ C_2(m,L) > 0.
\end{equation}

\vspace{4mm}

 \ A very popular example of these spaces  forms the so - called subgaussian space  $ \ \Sub = \Sub(\Omega); \ $ it consists on the
 subgaussian random variables, for which $ \  \psi(p) = \psi_2(p):= \sqrt{p}: \ $

 \begin{equation} \label{def sub}
 ||\tau|| \Sub = ||\tau||G\psi_2 \stackrel{def}{=} \sup_{p \ge 1} \left[   \  \frac{||\tau||_{p, \Omega}}{\sqrt{p}}  \ \right].
 \end{equation}
\ The r.v. $ \ \tau \ $ belongs to the subgaussian space $ \ \Sub(\Omega) \ $ iff

\begin{equation} \label{tail sub}
\exists C > 0 \ \Rightarrow \  T_{\tau}(u)  \le \exp(- C u^2), \ u \ge 0.
\end{equation}

\vspace{4mm}

 \ {\bf Example 4.1.}  \ Introduce the following $ \ G\Psi \ $ function

\begin{equation} \label{gamma example}
\nu[\gamma](p) = \nu(p) := \exp( 0.5 \gamma \ p ), \ p \ge 1, \ \gamma = \const  > 0.
\end{equation}

 \ If the r.v. $ \ \zeta \ $ belongs to the space $ \ G\nu[\gamma] \ $ and has therein an unit norm: $ \ ||\zeta||G\nu[\gamma] = 1, \ $
then

\begin{equation} \label{tail gamma}
T_{\zeta}(t) \le \exp \left( \  - 0.5 \ \gamma^{-1} \ (\ln^2 t) \ \right), \ t \ge e.
\end{equation}

\vspace{3mm}

 \ Conversely, let the estimation (\ref{tail gamma})  holds true for some r.v. $ \ \zeta; \ $ then this r.v. $ \ \zeta \ $ belongs to
 the Grand Lebesgue Space $ \ G \nu: \   ||\zeta||G\nu[\gamma]  \le C_1(\gamma) < \infty. \ $  \par

\vspace{4mm}

 \ {\bf Remark 4.1.} As a rule, on the the r.v.  $ \ \tau \ $ from the spaces $ \ G\psi_{m,L} \ $ is imposed the condition
 of {\it centering:} $ \ {\bf E} \tau = 0. \ $ \par

 \vspace{4mm}

 \ Another examples. Suppose that the r.v. $ \ \tau \ $ be such that

$$
T_{\tau}(t) \le T^{\beta,\gamma,L}(t), \ \beta > 1, \ \gamma > -1, \ L = L(t),
$$
where

$$
T^{\beta,\gamma,L}(t) \stackrel{def}{=} t^{-\beta} \ (\ln t)^{\gamma} \ L(\ln t), \ t \ge e
$$
and $ \ L = L(t), \ t \ge e \ $ be as before slowly varying at infinity positive continuous function.
It is known  \ \cite{KosOs equivalence} \ that as $ \ p \in [1,\beta) \ $

\begin{equation} \label{beta gamma}
\psi_{\tau}(p) = ||\tau||_p \le C_1(\beta,\gamma,L) \ (\beta - p)^{-(\gamma + 1)/\beta} \ L^{1/\beta}(1/(\beta - p)),
\end{equation}
and conversely, if the relation (\ref{beta gamma}) there holds, then

$$
T_{\tau}(t) \le  C_7(\beta,\gamma, L) \ T^{\beta,\gamma+1,L}(t).
$$

 \ Herewith both this estimations are unimprovable. \par

\vspace{5mm}

 \hspace{3mm} Let us return to the formulated above in this section problem. Indeed, we assume
 that the r.v. $ \ X \ $ belongs to the certain Grand Lebesgue Space (GLS) $ \ G\psi = G\psi(a,b): \ $

$$
||X||G\psi =  ||X||G\psi(a,b) < \infty; \hspace{3mm} 1 \le a < b \le \infty.
$$

 \ Of course, this generating function $ \ \psi(\cdot) \ $ may be choosed as a natural for the r.v.
 $ \ X: \ \psi_0(p) := |X|_p, \ $ if it is finite.\par
 \ Let $ \ \Delta = \const \in [a,b]; \ $ we introduce a new generating function

\begin{equation} \label{psi Delta}
\psi_{\Delta, \beta}(p) = \psi_{\Delta, \ \beta}[\psi](p)  \stackrel{def}{=} \frac{p}{p - \Delta} \ \beta^p \ \psi(p), \ \Delta < p \le b, \ \beta > 1.
\end{equation}
so that $ \ \psi_{\Delta,\beta} (\cdot) \in \Psi(\Delta,b). \ $ \par

\vspace{4mm}

 \ {\bf Proposition 4.1.}  One has in these notations, definitions  and under  our condition (\ref{cond Delta})

\vspace{3mm}

\begin{equation} \label{Proposit}
||\xi||G\psi_{\Delta,\beta } \le C \ ||X||G\psi,
\end{equation}

\vspace{3mm}

 with correspondent tail estimation \ (\ref{tail estim}).  \ Herewith the last estimation (\ref{Proposit}) is in general case essentially
 non -improvable.

\vspace{4mm}

 \ {\bf Proof.} One can take without loss of generality $ \ ||X||G\psi = 1; \ $ then

$$
\forall p \in (\Delta,b)  \ \Rightarrow |X|_p \le \psi(p).
$$

 \ We apply the estimation (\ref{beta Delta}) for these values $ \ p: \ $

$$
|\xi|_p \le C \ \frac{p}{ p - \Delta} \ \beta^p \ \psi(p) = C \ \psi_{\Delta,\beta}(p),
$$
or equally by means of the direct definition of the norm in the Grand Lebesgue Space $ \ G\psi_{\Delta} \ $

$$
||\xi||G\psi_{\Delta,\beta} \le C = C \ ||X||G\psi.
$$

 \vspace{3mm}

 \ The non - improvability of obtained estimate may be ground as before, as in the proposition 3.1,
 by means of considering at the example ones analogously as in the example 4.1:

 \begin{equation} \label{example  gamma}
T_{\xi}(t) = \exp \left( \  - 0.5 \ \gamma^{-1} \ (\ln^2 t) \ \right), \ t \ge e,
\end{equation}
where $ \ \gamma = 2 \ \ln \beta, \  \beta = \const > 1; \ $ and $ \ X =  1. \ $ \par

\vspace{3mm}

 \ In detail, it is easily to verify that the inequality (\ref{cond Delta}) for these variables
  is satisfied for all the values $ \ \beta > 1. \ $ It remains to note  as before that

$$
\sup_{\Delta > 1} \ \frac{||\xi||G\psi_{\Delta,\beta }}{ ||X||G\psi} = 1.
$$

 \ See for details the relation (\ref{beta Delta}).\par

\vspace{5mm}

\section{Multivariate case.}

\vspace{5mm}

 \hspace{3mm} We extend obtained results on the multidimensional case. Denotations and restrictions:

 $$
  d = \dim t = \dim{\vec{t}}  = \dim \xi  = \dim{\vec{\xi}} = 2,3,\ldots; \ \beta, C, \ \Delta = \const > 0;
 $$

$$
\vec{t} = \{ \ t_1, t_2, \ldots, t_d \ \}, \hspace{3mm} \vec{\xi} = \{ \ \xi_1, \xi_2, \ldots, \xi_d \ \};
$$

$$
\forall j = 1,2,\ldots, d \ \Rightarrow t_j, \ \xi_j \ge 0; \hspace{3mm} \vec{\xi} > \vec{t} \ \Leftrightarrow \forall j \ \xi_j > t_j;
$$

$$
|\vec{t}| \stackrel{def}{=} \sqrt{ \ \sum_{j=1}^d t_j^2 \ },  \hspace{3mm} |\vec{\xi}|_p \stackrel{def}{=} \left[ \ \sum_{j=1}^d |\xi_j|_p^p  \ \right]^{1/p}, \ p \ge 1.
$$

\vspace{4mm}

 \hspace{3mm} {\bf  Given:}

\begin{equation} \label{d datum}
| \ \vec{t} \ |^{\Delta} \ {\bf P} (\vec{\xi} >  \beta \vec{t})   \le C \ {\bf E}[ \ X \ I(\vec{\xi} > \vec{t}) \ ].
\end{equation}

\vspace{4mm}

\ {\bf Proposition 5.1.} We state for the values $ \ p > \max(1,\Delta) \ $

\vspace{4mm}

\begin{equation} \label{d dim est}
|\vec{\xi}|_p \le C \ \frac{p}{p - \Delta} \ d^{1/p} \ \beta^p \ |X|_p.
\end{equation}

\vspace{4mm}

\ {\bf Proof.}  We have  for the values $ \ j = 1,2,\ldots,d \ $

$$
t_j \ {\bf P} (\xi_j > \beta \ t_j)  \le C {\bf E} [ \ X I(\xi_j > t_j) \ ].
$$
 \ It follows from the one - dimensional estimates

$$
|\xi_j|_p \le C \ \frac{p}{p - \Delta} \  \beta^p \ |X|_p.
$$

 \ Further,

$$
|\vec{\xi}|_p = \left[ \ \sum_{j=1}^d |\xi_j|_p^p \ \right]^{1/p}  = \left[ \ \sum_{j=1}^d |\xi_j|_p^p  \ \right]^{1/p} \le
$$

$$
\left[ \ d \ \left( \ \frac{C p}{p - \Delta} \ \right)^p  \ \beta^{p^2} \ |X|_p^p \ \right]^{1/p}  =
C \ d^{1/p} \ \frac{p}{p - \Delta} \ \beta^p \ |X|_p,
$$

\vspace{3mm}

 \ Q.E.D. \par

\vspace{4mm}

 \ {\bf Remark 5.1.} One can use instead the classical Euclidean norm for the vector $ \ t: \ |t| \ $ arbitrary another one,
 for instance, $ \ l_s \ $ one:

$$
|t|_s := (\sum_{j=1}^d \ |t_j|^s)^{1/s},  \ s = \const \ge 1.
$$

 \vspace{5mm}

\section{Concluding remarks.}

 \vspace{5mm}

 \hspace{3mm} It is interest in our opinion to apply obtained in this preprint estimates, in particular, in the
 martingale theory.\par

\vspace{6mm}

\vspace{0.5cm} \emph{Acknowledgement.} {\footnotesize The first
author has been partially supported by the Gruppo Nazionale per
l'Analisi Matematica, la Probabilit\`a e le loro Applicazioni
(GNAMPA) of the Istituto Nazionale di Alta Matematica (INdAM) and by
Universit\`a degli Studi di Napoli Parthenope through the project
\lq\lq sostegno alla Ricerca individuale\rq\rq .\par

\end{document}